\documentclass[12pt]{article}
\usepackage{graphicx,amsmath,amsthm,amssymb,dsfont, mathrsfs,MnSymbol}
\usepackage{setspace,bbm}
\def\ZZ{\mathbb Z}

\def\QQ{\mathbb Q}

\def\DD{\mathbb D}

\def\VV{\mathbb V}
\def\UU{\mathbb U}

\def\cH{\mathcal H}

\def\ff{\mathfrak f}

\def\fm{\mathfrak m}

\def\mm{\mathbbm m}

\def\fd{\mathfrak d}

\def\fO{\mathcal O}
\def\cP{\mathcal P}

\def\cD{\mathcal D}

\def\sD{\mathscr D}
\def\cN{\mathcal N}

\def\fT{\mathfrak T}

\def\bD{{\bf D}}

\def\bx{{\bf x}}
\def\by{{\bf y}}

\def\bn{{\bf n}}
\def\bm{{\bf m}}

\def\bb{{\bf b}}

\def\br{{\bf r}}
\def\bs{{\bf s}}

\def\b1{{\bf 1}}
\def\d1{\mathds{ 1}}
\def\mod{{\; \rm mod \;}}
\def\ad{{\rm and}}

\def\with{{\rm with}}
\def\where{{\rm where}}
\def\for{{\rm for}}

\def\ord{{\rm ord}}

\def\qed{ \ \vrule width.2cm height.2cm depth0cm\smallskip}

\begin{document}


\title{On the upper bound of the $L_2$-discrepancy  of Halton's sequence}
\author{Mordechay B. Levin}

\date{}

\maketitle

\begin{abstract}
Let $ (H(n))_{n \geq 0} $ be a $2-$dimensional Halton's sequence. Let \\
 $D_{2} ( (H(n))_{n=0}^{N-1}) $ be the $L_2$-discrepancy of  $ (H_n)_{n=0}^{N-1} $.
It is known that \\ $\limsup_{N \to  \infty }
 (\log N)^{-1}  D_{2}  (  H(n) )_{n=0}^{N-1} >0$.
 In this paper, we prove that
 $$
           D_{2} ((  H(n) )_{n=0}^{N-1})  =O(  \log N) \quad {\rm for} \; \; N \to  \infty ,
 $$
i.e., we found the smallest possible order of magnitude of $L_2$-discrepancy  of a 2-dimensional Halton's sequence.
The main tool is  the theorem on linear forms in the $p$-adic logarithm.
\end{abstract}
%
{\bf 1. \;\; Introduction }\\
%
 Let $\cP_{N} = (\beta_{n})_{n = 0}^{N-1}$ be an $N$-element point set in the
 $s$-dimensional unit cube $[0,1)^s$. The local {\bf discrepancy} function of $\cP_{N}$ is defined as
\begin{equation}\label{In1}
   D(\bx, \cP_{N}  )= \sum\nolimits_{n=0}^{N-1}   \b1_{B_{\bx}}(\beta_{n}) - N x_1 \cdots x_s,
\end{equation}
where $B_{\bx}=[0,x_1) \times \cdots \times [0,x_s) $,
 $\b1_{B_{\bx}}(\by) =1, \; {\rm if} \;\by  \in B_{\bx}$,
and $   \b1_{B_{\bx}}(\by) =0,$  if $
\by   \notin B_{\bx}$.
We define the $\emph{L}_q$  discrepancy of an
$N$-point set $\cP_{N}$ as
\begin{equation} \label{In2}
       D_{\infty} (\cP_{N}) =
    \sup_{ 0<x_1, \ldots , x_s \leq 1} \; |
    D(\bx,\cP_{N}) |,\qquad  D_{q}(\cP_{N})=\left\| D(\bx,\cP_{N}) \right\|_{q}  ,
\end{equation}
\begin{equation*} 
     \left\| f(\bx) \right\|_{q}= \Big( \int_{[0,1]^s}  |f(\bx)|^q d \bx \Big)^{1/q}.
\end{equation*}

 A sequence $(\beta_n)_{n\geq 0}$ is of {\it
low discrepancy} (abbr. l.d.s.) if  $D_{\infty}
((\beta_n)_{n=0}^{N-1})=O(\log^s N) $.
A point set $\cP_{N}$ is of
 low discrepancy (abbr.
l.d.p.s.) if $ D_{\infty} (\cP_{N})=O(\log^{s-1}
N)$.

Let $p_1,\ldots , p_s \geq 2 $ be pairwise coprime integers,
 \begin{equation}\label{In5}
 n=\sum_{j\geq 1}e_{i,j}(n) p_i^{j-1},\;  e_{i,j}(n) \in \{0,1, \ldots
 ,p_i-1\}, \;  {\rm and}    \; \phi_i(n)= \sum_{j\geq 1}e_{i,j}(n) p_i^{-j}.
 \end{equation}
%
%
Van der Corput    proved that $ (\phi_1(n))_{n\geq 0}$ is the $1-$dimensional l.d.s.\\
The first example of multidimensional l.d.s. was proposed by Halton
\begin{equation}  \label{In6}
  H_s(n)= (\phi_{1}(n),\ldots ,\phi_{s}(n)), \quad n=0,1,2,... \;.
\end{equation}
The first example of multidimensional l.d.p.s. was obtained by Hammersley
${\cH}_{s+1,N}=
	 (H_s(n), n/N)_{n=0}^{N-1}$.
For other examples of  l.d.s. and  l.d.p.s., see, e.g., [BC],  [Ni].
It is known that
\begin{equation} \label{In7}
\emph{D}_{s,q}(\cP_N)_{n=0}^{N-1})>C_{s,q}(\log N)^{\frac{s-1}{2}}, \;\;\;\;  \overline{\lim }{\emph{D}_{s,q}((\beta_n)_{n=0}^{N-1})(\log
N)^{-s/2}}>0
\end{equation}
for all point sets $\cP_N$ and all sequences
$(\beta_n)_{n \geq 0}$ with some $C_{s,q} >0$ (see  Roth for $q=2$, Schmidt  for $q>1$ [BeCh], and Proinov [Pr]).

 A sequence $(\beta_n)_{n\geq 0}$ is of  $L_q${\it -low discrepancy} if  $ D_{q}
((\beta_n)_{n=0}^{N-1})=O(\log^{s/2} N) $.
A point set $\cP_{N}$ is of  $L_q$-{\it
low discrepancy} if $ D_{q} (\cP_{N})=O( \log^{(s-1)/2}
N) $.

The existence of $L_2$-l.d.p.s. was proved by Roth [Ro]. For the proof, Roth used the averaging of
 Hammersley's point set.
The existence of $L_q$-l.d.p.s. for $q>1$ was proved
 by Chen   [Ch].
The first explicit construction of $L_q$-l.d.p.s. was obtained by Chen and Skriganov for $q=2$ and by Skriganov  for $q>1$ (see [ChSk], [Sk]).
The next explicit construction of $L_q$-l.d.p.s. was proposed by Dick and Pillichshammer
 (see [Di], [DP],  [Ma]).
The first explicit construction of $L_q$-l.d.s. were obtained by Dick, Hinrichs, Markhasin and Pillichshammer [DHMP]. All these explicit constructions were obtained by using $(t,m,s)$ nets.
In [Le], we obtained a similar result for Halton's sequence. The proof of [Le] is very technical.
So, to explain the approach of [Le], in this paper we will consider only the case  $s = 2, \; q =2$ :\\ \\
{\bf Theorem.} {\it Let $p_1,p_2$ be primes, $p_1 \neq p_2$ and let $Q \geq 0, N \geq 2$ be integers. Then}
\begin{equation} \nonumber
  D_{2} ((  H_{2}(k) )_{k=Q}^{Q+N-1})  \leq C \log N   ,
\end{equation}
where the constant $C$ is independent of $Q$ and $N$. \\

Now we describe the structure of the paper. In Lemma 1 - Lemma 3, we get a  simple estimate of Fourier's series of truncated discrepancy function of Halton's sequence  $\cP_{N} =( H_{2}(k) )_{k=Q}^{Q+N-1}$.
Next, we divide  the box $ [1, n] ^ 4$  and $  D^2 (\bx,\cP_{N})$ to four partitions :
\begin{equation}\nonumber
  D^2 (\bx,\cP_{N}) = \bigcup_{\lambda_i=0,1,\;i=1,2}
    \bD_{\lambda_1, \lambda_2,N}, \qquad
  D^2_{2} (\cP_{N})\leq \sum_{\lambda_i=0,1,\;i=1,2}
    \left\| \bD_{\lambda_1, \lambda_2,N} \right\|_{1},
\end{equation}
$ [1, n] ^ 4 =\bigcup_{\lambda_i=0,1,\;i=1,2} \UU_{ \lambda_1, \lambda_2 } $, $n=[\log_2 N] +1$,      and we prove
\begin{equation}\nonumber
      \left\| \bD_{\lambda_1, \lambda_2,N} \right\|_{1}  \ll
\sum_{ (\br_1,\br_2) \in \UU_{ \lambda_1, \lambda_2 }} \tilde{\DD}_{\lambda_1, \lambda_2,\br_1, \br_2},
  \quad  \br_j =(r_{1,j},r_{2,j})\in [1,n]^2, \; j=1,2,
\end{equation}
where $\tilde{\DD}_{\lambda_1, \lambda_2,\br_1, \br_2} $ is the interior sum in Lemma 4. In Lemma 5, we  minimise the number of terms in the expression of Lemma~4.
The main tool of the proof is the  Yu theorem on linear forms in $p$-adic logarithm (see [Yu], [Bu]).
In Lemma 6 and Lemma 7, we apply this theorem  to get the estimate $ \left\| \bD_{\lambda_1, \lambda_2,N} \right\|_{1} \ll n^2$.
\\ \\
{\bf 2. \;\; Beginning of the proof of the Theorem.}
 We will use notation $A \ll B$  equal to $A = O(B)$.\\ \\
{\bf Lemma 1.} {\it For each $i \in \{1,2\}$, suppose that $y_i \in [0,1)$ and $k_i \in \ZZ$ satisfy             \begin{equation} \nonumber
   y_i= \sum_{1 \leq j \leq s_i} y_{i,j} p_i^{-j} \quad \ad \quad
       k_i= \sum_{1 \leq j \leq s_i} y_{i,j} p_i^{j-1} ,
\end{equation}
where $y_{i,j} \in \{ 0,1,...,p_i-1  \}$ for every $ 1 \leq j \leq s_i $,\\
\quad $(i)\;$  For each $i \in \{1,2\}$, we have
\begin{equation} \nonumber
   \phi_i(k) \in [y_i, y_i + p_i^{-s_i}) \quad {\rm  if  \; and \; only \; if} \quad k \equiv k_i \mod p_i^{s_i}.
\end{equation}
\quad $(ii)\;$  We have
\begin{equation} \nonumber
  H_2(k) =( \phi_1(k), \phi_2(k)) \in [y_1, y_1 + p_1^{-s_1}) \times [y_2, y_2 + p_2^{-s_2})
\end{equation}
\begin{equation} \nonumber
   \quad {\rm  if  \; and \; only \; if} \quad \quad k \equiv
   p_2^{s_2} M_{1,\bs} k_1 +  p_1^{s_1} M_{2,\bs} k_2 \mod p_1^{s_1}p_2^{s_2},
\end{equation}
where for $\bs =(s_1,s_2)$, the integers $M_{1,\bs}$ and  $M_{2,\bs}$ satisfy the congruences}
\begin{equation} \label{Beg-3} 
  p_2^{s_2} M_{1,\bs} \equiv 1 \mod p_1^{s_1}, \;\; \ad \;\;
p_1^{s_1} M_{2,\bs} \equiv 1 \mod p_2^{s_2}, \; M_{i,\bs} \in [0, p_i^{s_i}).
\end{equation}
{\bf Proof.} $(i)$ follows from the definition of $\phi_i(k)$ given by \eqref{In5},
and $(ii)$ is then a special case of the Chinese Remainder Theorem. \qed \\

For $i \in \{1,2\}$, we can write $x_i \in [0,1)$ in base $p_i$ expansion in the form
 $x_i =0.x_{i,1}x_{i,2}...=\sum_{j \geq 1}  x_{i,j} p_{i}^{-j}$, with $x_{i,j} \in \{0,1,...,p_{i}-1  \}$.
For any positive integer $r_i$, we write
\begin{equation}  \nonumber
        [x_i]_{r_i} =\sum_{1 \leq j \leq r_i} x_{i,j} p_{i}^{-j}
\end{equation}
to denote the expansion of $x_i$ truncated after $r_i$ digits. Furthermore, the truncation
 $[\bx]_{\br}$, where $\br =(r_1,r_2)$, is defined coordinatewise by writing $[\bx]_{\br}=
( [x_1]_{r_1},  [x_2]_{r_2})$ with  $\bx = (x_1,  x_2)$.
  A special case of Lemma 1 is that for $\br =(r_1,r_2)$, and $P_{\br} =  p_1^{r_1}p_2^{r_2} $, we have
\begin{multline}\label{Beg-6} 
  [H_2(k)]_{\br} = [\bx]_{\br}  \;\;\; {\rm  if  \; and \; only \; if}    \;\;\;  k \equiv   X_{\br}  \; {\rm mod} \; P_{\br}, \quad \where\\
 \quad  X_{\br} \equiv    p_2^{r_2} M_{1,\br}
    X_{1,r_1} +   p_1^{r_1} M_{2,\br}
    X_{2,r_2} \;\mod \; P_{\br}, \;\;\;\; X_{\br} \in [0,P_{\br}), \\
    X_{i,r_i}=
  \sum_{1 \leq j \leq r_i}  x_{i,j} p_i^{j-1}, \;  \;
	 \quad X_{i,r_i} \in [0,p_i^{r_i}), \;\; 1 \leq i \leq 2.
\end{multline}

Let
\begin{equation} \nonumber
   \delta_M(a) =   \begin{cases}
    1,  & \; {\rm if}  \;  a \equiv 0  \mod M,\\
    0, &{\rm otherwise},
  \end{cases}  \qquad \ad \qquad
   \Delta(\fT) =   \begin{cases}
    1,  & \; {\rm if}  \;  \fT  \;{\rm is \;true},\\
    0, &{\rm otherwise}
  \end{cases}.
\end{equation}

Let $n =[\log_2 N]+1$ and let $\bn =(n,n)$. We shall make an approximation to the function
$ D(\bx, (H_2(k))_{k=Q}^{Q+N-1}  )  $ by the quantity
\begin{equation}  \label{Beg-16}
  \sD(N) :=  D([\bx]_{\bn}, (H_2(k))_{k=Q}^{Q+N-1}  ) .
\end{equation}

Repeating the proof from [Ni, p.29-31], we obtain from \eqref{In1} the following lemma :\\ \\
{\bf Lemma 2.} {\it We have
\begin{equation}  \label{Beg-22}
     \sD(N)  =\sum_{r_1,r_2 =1}^n   \cD_{\br,N}=  D(\bx, (H_2(k))_{k=Q}^{Q+N-1}  ) + 2\epsilon ,  \;\; |\epsilon| \leq 1, \quad \where
\end{equation}
\begin{equation}  \label{Beg-18a}
  \cD_{\br,N} = \sum_{b_1=0}^{x_{1,r_1}-1}  \sum_{b_2=0}^{x_{2,r_2} -1}  \sum_{k=Q}^{Q+N-1} (\delta_{P_{\br}} (k - X_{\br,\bb})-1/P_{\br} ),
\end{equation}
\begin{equation}  \label{Beg-18}
  X_{\br,\bb} \equiv \sum_{i=1}^2  M_{i,\br}
   P_{\br}  p_i^{-r_i}
  \Big( \sum_{1 \leq j < r_i}  x_{i,j} p_i^{j-1}  + b_i p_i^{r_i-1} \Big) \mod  P_{\br} , \quad
	X_{\br,\bb} \in [0,P_{\br}).
\end{equation}
Furthermore, $ | \cD_{\br,N}| \leq x_{1,r_1} x_{2,r_2}  <p_1 p_2$ for every $\br=(r_1,r_2)$ and every
natural number $N$. Here we use the convention that the sum \eqref{Beg-18} is empty if $x_{1,r_1} =0$ or
$x_{2,r_2} =0$. }\\

Let $[y]$ be the integer part of $y$,
\begin{equation} \label{Beg-1a}
   I_M=[-[(M-1)/2],[M/2]] \cap \ZZ \quad \ad \quad  I_M^{*}=I_M \setminus \{ 0 \}.
\end{equation}
Then  $I_M$ is a complete set of residues $\mod M$, $M \geq 1$. It is well known that
\begin{equation} \label{Beg-1}
   \delta_M(a) =   \frac{1}{M} \sum_{k\in I_M} e\Big(\frac{ak}{M}\Big),  \qquad \where \quad e(x) =exp(2\pi i x),
\end{equation}
and that
\begin{equation} \label{Beg-2a}
  \Big|\frac{1}{M} \sum_{k=0}^{M-1} e\big(k \alpha) \Big| \leq \min\Big(1, \frac{1}{2M\llangle \alpha \rrangle} \Big) ,\quad \;\; \where \;
	\llangle \alpha \rrangle = \min( \{\alpha\}, 1 - \{\alpha\})
\end{equation}
see, e.g., [Ko, Lemma 2, p. 2] and [Ko, Lemma 1, p. 1].

Let $\bar{m} = \max(1, |m|) \leq M/2 $.     From [Ko, p. 2], we obtain for $R \leq M$ :
\begin{equation} \label{Beg-2}
  \Big|\frac{1}{M} \sum_{k=0}^{R-1} e\big(\frac{m k}{M}\big) \Big| \leq \min\Big(1, \Big|\frac{e(m R/M) -1}{M
	(e(m /M)-1)} \Big|\Big) \leq \frac{|\sin(\pi m R/M)|}{\bar{m}}  \leq \frac{1}{\bar{m}}.
\end{equation}\\
{\bf Lemma 3.} {\it With the notations as above, we have}
\begin{eqnarray}
&   \cD_{\br,N} =
  \sum_{m \in I_{P_{\br}}^{*}}   \varphi_{ \br,Q,N,m} \; \psi_{ \br}(m,\bx)\;
e\Big( \frac{-m}{P_{\br}} X_{\br} \Big),  \label{Le2-1} \\
&where \; X_{\br} \; is \; given \;by \;\eqref{Beg-6} \;and \qquad \qquad\qquad \qquad\qquad \qquad\qquad \qquad\qquad   \nonumber\\
& \varphi_{ \br,Q,N,m} =\frac{1}{P_{\br}}
 \sum_{k=Q}^{Q+N-1} e\Big( \frac{m k}{P_{\br}} \Big) ,  \quad  satisfies \quad
  |\varphi_{ \br,Q,N,m}| \leq \frac{1}{\bar{m}}, \label{Le2-51}\\
&and \qquad \qquad\qquad \qquad\qquad \qquad\qquad  \qquad\qquad \qquad\qquad \qquad
\qquad \qquad \nonumber\\
&        \psi_{ \br}(m,\bx)  = \prod_{i=1}^2 \tilde{\psi}(i,m, x_{i, r_i}) ,\quad satisfies   \quad |\psi_{ \br}(m,\bx)| \leq p_1 p_2, \label{Le2-52} \\
&where, \;for \;each \;i \in \{1,2\},  \qquad \qquad\qquad \qquad\qquad \qquad\qquad \qquad\qquad \nonumber\\
& 			\tilde{\psi}(i,m, x_{i, r_i}) =
\sum_{b_i=0}^{x_{i,r_i} -1}
e\Big(\frac{-m M_{i,\br}(b_i-x_{i, r_i})}{p_i}\Big) \label{Le2-53} .
\end{eqnarray} \\
{\bf Proof.}
Combining \eqref{Beg-18a} and \eqref{Beg-1}, we have
\begin{multline*}
 \cD_{\br,N}=    \sum_{b_1=0}^{x_{1,r_1} -1} \; \sum_{b_2=0}^{x_{2,r_2} -1} \sum_{k=Q}^{Q+N-1} \frac{1}{P_{\br}} \sum_{m \in I_{P_{\br}}^{*}} e\Big( \frac{m}{P_{\br}} (k - X_{\br,\bb})\Big) \\
=    \sum_{m \in I_{P_{\br}}^{*}}
 \Bigg( \frac{1}{P_{\br}}
 \sum_{k=Q}^{Q+N-1} e\Big( \frac{m k}{P_{\br}} \Big)   \Bigg)
\; \sum_{b_1=0}^{x_{1,r_1} -1}   \sum_{b_2=0}^{x_{2,r_2} -1} e\Big( \frac{-m}{P_{\br}} X_{\br,\bb} \Big) \\
=    \sum_{m \in I_{P_{\br}}^{*}} \varphi_{ \br,Q,N,m}
\; \sum_{b_1=0}^{x_{1,r_1} -1}   \sum_{b_2=0}^{x_{2,r_2} -1} e\Big( \frac{-m}{P_{\br}} X_{\br,\bb} \Big)               .
\end{multline*}
Comparing \eqref{Beg-6} and \eqref{Beg-18}, we see that
\begin{equation*}  
  X_{\br,\bb} \equiv X_{\br} +    \frac{P_{\br} M_{1,\br} }{p_1}   (b_1 - x_{1,r_1}) +
   \frac{P_{\br} M_{2,\br} }{p_2}   (b_2 - x_{2,r_2})
   \mod  P_{\br}.
\end{equation*}
It follows that
\begin{multline*}
 \sum_{b_1=0}^{x_{1,r_1} -1}  \sum_{b_2=0}^{x_{2,r_2} -1} e\Big( \frac{-m}{P_{\br}} X_{\br,\bb} \Big)  =
\sum_{b_1=0}^{x_{1,r_1} -1}   \sum_{b_2=0}^{x_{2,r_2} -1}
e\Big( \frac{-m}{P_{\br}} \big( X_{\br} +    \frac{P_{\br} M_{1,\br} }{p_1}   (b_1 - x_{1,r_1}) \\
 +
   \frac{P_{\br} M_{2,\br} }{p_2}   (b_2 - x_{2,r_2})\big) \Big) =  \psi_{ \br}(m,\bx) \;
e\Big( \frac{-m}{P_{\br}} X_{\br} \Big).
\end{multline*}
establishing \eqref{Le2-1}. Finally, the estimate \eqref{Le2-51} follows immediately from \eqref{Beg-2},
 while the estimate in \eqref{Le2-52} follows on observing that each summand  in \eqref{Le2-53} is bounded by 1. Hence Lemma 3 is proved. \qed \\
%

%
Let
\begin{multline}   \label{Beg-26}
  V =p_1^4 p_2^4 [\log_2^4 n], \;\; U= \{ (r_1,r_2) \in [1,n]^2 \;  | \;  \max(r_1,r_2) >V\} \quad  \ad \; {\rm write}\\
             \tilde{\sD}(N)=  \sum_{\br \in U}    \cD_{\br,N}.
             \quad\quad\quad\quad\quad\quad\quad\quad \quad\quad
\end{multline}
Here we cut out a small box using the parameter $V$ for the purposes of Lemma 5.
By Lemma 2, we have $|\cD_{\br,N}| < p_1p_2$.
From \eqref{Beg-22} and  \eqref{Beg-26}, we derive
\begin{equation}  \nonumber 
| \sD(N) -  \tilde{\sD}(N) | \leq  p_1 p_2V^2 \leq p_1^5 p_2^5 \log_2^{8} n.
\end{equation}
Using  \eqref{Beg-22}, we obtain
\begin{equation}  \label{Beg-28}
 D(\bx, (H(k))_{k=Q}^{Q+N-1}  )   =   \tilde{\sD}(N)    + O(\log^{8} n).
\end{equation}

Let
\begin{multline}   \label{Beg-40}
 \VV =[\log_2^3 n],\;\;
  \br_j =(r_{1,j},r_{2,j}), \;\; r^{+}_i =\max_j  r_{i,j}, \; \;r^{-}_i =\min_j  r_{i,j}, \\
  \UU_{\lambda_1, \lambda_2}= \{ (\br_1,\br_2) \in U^2 \; | \;\;  \Delta(r^{+}_i - r^{-}_i > \VV )  =\lambda_i, \;\; i=1,2 \},\\
%
  \bD_{\lambda_1, \lambda_2,N}= \sum_{(\br_1,\br_2) \in \UU_{\lambda_1, \lambda_2}}   \cD_{\br_1,N} \cD_{\br_2,N}, \quad \lambda_i \in \{0,1\}, \;\; i=1,2.
\end{multline}
It is easy to verify that $U^2 = \bigcup_{\lambda_i \in \{0,1\}, i=1,2}  \UU_{\lambda_1, \lambda_2}$,
 $\UU_{\lambda_1, \lambda_2} \bigcap \UU_{\lambda_1^{'}, \lambda_2^{'}} =\emptyset $ for $(\lambda_1, \lambda_2) \neq (\lambda_1^{'}, \lambda_2^{'})$.
By   \eqref{Beg-26}, we have
\begin{equation}   \label{Beg-44}
\tilde{\sD}(N)^2  = \sum_{\substack{\lambda_i \in \{ 0,1 \} \\ i=1,2 }}  \bD_{\lambda_1, \lambda_2,N} \quad
    \ad \quad
   \left\| \tilde{\sD}(N) \right\|_{2}^2   \leq   \sum_{\substack{\lambda_i \in \{ 0,1 \} \\ i=1,2 }}
     \Big| \int_{[0,1]^2}  \bD_{\lambda_1, \lambda_2,N}    d \bx \Big| .
\end{equation}
 Here we divided the set $ U^2$ into four partitions using the parameter $\VV$ for the purposes of Lemma 6 and Lemma7.
From   \eqref{Beg-28} and \eqref{Beg-44}, we get that in order to prove the Theorem it is enough to verify that
\begin{equation}   \label{Beg-90}
  \Big| \int_{[0,1]^2}  \bD_{\lambda_1, \lambda_2,N}    d \bx \Big|     \ll  n^2, \quad \for \quad \lambda_i \in \{ 0,1 \}, \; i=1,2.
\end{equation}
Let $\br_j=(r_{1,j},r_{2,j})$,  $j \in [1,2]$,  $k_{i,j}$  satisfy the conditions
\begin{equation}\label{Le3-17}
               r_{i,k_{i,2}} = r^{+}_i =\max_j  r_{i,j}, \quad  r_{i,k_{i,1}} = r^{-}_i =\min_j  r_{i,j} \quad \for \quad r^{+}_i > r^{-}_i,
\end{equation}
 $(k_{i,1},k_{i,2})=(1,2) $ for  $r^{+}_i = r^{-}_i$, and let $k_{i,0} = r_{i,0} =0$
  $(i \in [1,2] )$.\\
%
%
Let
\begin{equation} \label{Le3-9a}    
\hat{m}_{i}  \equiv  - \sum_{j=1}^{2}   m_j M_{i,\br_j}
   p_i^{r^{+}_i -r_{i,j} } \mod p_i^{r^{+}_i}, \; \quad
       \hat{m}_{i} \in [0,p_i^{r^{+}_i})    .
\end{equation}
Let $\mm_{i,j}$  be the unique
integer satisfying the following conditions
 (see \eqref{Beg-1a}) :
\begin{equation} \label{Le3-20}
\hat{m}_{i}  \equiv \sum_{j=1}^{2}  \mm_{i,j} p_i^{ r^{+}_i -  r_{i,j}  } \; \mod p_i^{r^{+}_i}  , \quad \with  \quad    \mm_{i,k_{i,j}} \in I_{p_i^{ r_{i,k_{i,j}}  -r_{i,k_{i,j-1}} }},
\end{equation}
 $ i,j=1,2, \; I_{1} = I_{p_i^0}=\{ 0 \} $. \\
Let $E_i(f(x_i))= \int_{[0,1]} f(x_i)d x_i  $.  It is easy to see that
\begin{equation} \label{Le3-9}
 E_i(f([x_i]_l))=  \frac{1}{p_i^l} \sum_{x_{i,1}=0}^{p_i-1}...\sum_{x_{i,l}=0}^{p_i-1} f(0.x_{i,1} ... x_{i,l}).
\end{equation} \\
{\bf Lemma 4.}  {\it With the notations as above, we get}
\begin{equation}  \label{Le4-0}
 \Big| \int_{[0,1]^2}  \bD_{\lambda_1, \lambda_2,N}   d \bx \Big|   \ll
\sum_{ (\br_1,\br_2) \in \UU_{ \lambda_1, \lambda_2 }} \;\;\;
\sum_{  m_j \in I^{*}_{P_{\br_j}}, \;  j=1,2} \;
\frac{ 1 }{\bar{m}_1  \bar{m}_{2}}
 \prod_{i,j=1}^{2} \frac{1}{\bar{\mm}_{i, k_{i,j}}} .
\end{equation}  \\
{\bf Proof.}
From \eqref{Beg-6} and \eqref{Le3-9a}, we have
\begin{multline*}   
 e\Big(-\sum_{j=1}^{2}  \frac{m_j }{P_{\br_j}}X_{\br_j}  \Big)
	  =
 e\Big(-  \sum_{i,j=1}^{2}  \frac{m_j M_{i,\br_j} X_{i,r_{i,j}}}{p_i^{r_{i,j}}}  \Big)
  	  \\
=   e\Big(- \sum_{i,j=1}^{2}   \frac{m_j M_{i,\br_j} p_i^{r^{+}_i -r_{i,j}}
 }{p_i^{r^{+}_i}}
  X_{i,r^{+}_i} \Big)
 =   e\Big(  \sum_{i=1}^2 \frac{\hat{m}_{i} X_{i,r^{+}_i}}{p_i^{r^{+}_i}}
	\Big).
\end{multline*}
By  \eqref{Le2-1} - \eqref{Beg-40}, we obtain that
$ \psi_{ \br}(m,\bx)  = \prod_{i=1}^2 \tilde{\psi}(i,m, x_{i, r_i})  $,  $ |\varphi_{ \br_j,Q,N,m_j}| \\ \leq 1/\bar{m}_j $ and
\begin{multline}  \nonumber 
 \Big| \int_{[0,1]^2}  \bD_{\lambda_1, \lambda_2,N}   d \bx \Big|
   \leq
  \sum_{ (\br_1,\br_2) \in \UU_{ \lambda_1, \lambda_2 }}
\sum_{\substack{m_j \in I^{*}_{P_{\br_j}} \\ j=1,2}}
 \Big| \int_{[0,1]^2} \prod_{j=1}^{2}
     \varphi_{ \br_j,Q,N,m_j} \; \psi_{ \br_j}(m_j,\bx)  \\
  \times e\Big( \frac{-m_j}{P_{\br_j}} X_{\br_j} \Big)     d \bx \Big|
\leq
\sum_{ (\br_1,\br_2) \in \UU_{ \lambda_1, \lambda_2 }}
\sum_{\substack{m_j \in I^{*}_{P_{\br_j}} \\ j=1,2}}
 \frac{1}{ \bar{m}_1  \bar{m}_{2}}  \Big| \int_{[0,1]^2}
   \prod_{j=1}^{2}  \psi_{ \br_j}(m_j,\bx)\;
e\Big( \frac{-m_j}{P_{\br_j}} X_{\br_j} \Big)    d \bx \Big|  \\
\leq
\sum_{ (\br_1,\br_2) \in \UU_{ \lambda_1, \lambda_2 }} \;
\sum_{ m_j \in I_{P_{\br}}^{*}, \; j=1,2} \;
 \frac{1}{ \bar{m}_1  \bar{m}_{2}} \; \prod_{j=1}^{2}\; \Big| \int_{[0,1]}  Z_i d x_i \Big| ,\\
  \; \with \;
 Z_i=  e\Big(   \sum_{i=1}^s \frac{\hat{m}_{i} X_{i,r^{+}_i}}{p_i^{r^{+}_i}}
	\Big)   \prod_{j=1}^{2}
     \tilde{\psi}\Big(i , m_{j}  , x_{i,r_{i,j}}  \Big).
\end{multline}
%
%

%
%
We get that in order to prove Lemma 4 it is enough to verify the following inequality
\begin{equation}  \label{Le2-20}
 \Big|  \int_{[0,1]}  Z_i d x_i \Big|  \leq
   \prod_{j=1}^{2} \frac{p_i^4}{\bar{\mm}_{i, k_{i,j}}}.
\end{equation}
In view of  \eqref{Beg-6} and \eqref{Le3-9}, we have
\begin{equation}  \label{Le2-01}
  \int_{[0,1]}  Z_i d x_i =E_i(Z_i) =\frac{1}{p_i^{r^{+}_i  }}   \sum_{\substack{x_{i,k} \in [0,p_i)\\k \in [1,r^{+}_i]}}       e\Big(  \frac{\hat{m}_{i} X_{i,r^{+}_i}}{p_i^{r^{+}_i}}
	\Big)   \prod_{j=1}^{2}      \tilde{\psi}(i , m_{j}  , x_{i,r_{i,j}}  ) ,
\end{equation}
%
%
where $ X_{i,r_i^{+}}=  \sum_{1 \leq j \leq r_i^{+}}  x_{i,j} p_i^{j-1} $.
According to \eqref{Le3-17}, we obtain $   r_{i,k_{i,2}} = r^{+}_i \\=\max_j  r_{i,j}$,  $r_{i,k_{i,1}} = r^{-}_i =\min_j  r_{i,j}, \; r_{i,0}=0$. We divide the interval $[1, r_i^{+}]$ into two parts
$[1,r_{i,k_{i,1}}]$ and $[r_{i,k_{i,1}} +1, r_{i,k_{i,2}}]$.
By    \eqref{Le2-01}, we get
\begin{multline*} \nonumber
  E_i(Z_i) =\frac{1}{p_i^{r^{+}_i  }}   \sum_{\substack{x_{i,k} \in [0,p_i)\\k \in [1,r^{+}_i]}}
       \prod_{j=1}^{2}
         e\Bigg(  \frac{\hat{m}_{i} \sum_{k=r_{i,k_{i,j-1}}+1}^{r_{i,k_{i,j}} } x_{i,k}p_i^{k-1}}{p_i^{r^{+}_i}}
	\Bigg)
	        \tilde{\psi}(i , m_{j} , x_{i,r_{i,j}}  ) \\
= \prod_{j=1}^{2} \frac{1}{p_i^{r_{i,k_{i,j}} -r_{i,k_{i,j-1}}}}
\sum_{\substack{x_{i,k} \in [0,p_i)\\
k \in [r_{i,k_{i,j-1} }+1,r_{i,k_{i,j}} ]}}
        e\Bigg(  \frac{\hat{m}_{i} \sum_{k=r_{i,k_{i,j-1}}+1}^{r_{i,k_{i,j}} } x_{i,k}p_i^{k-1}}{p_i^{r^{+}_i}}
	\Bigg)  \tilde{\psi}(i ,m_j , x_{i,r_{i,j}} ).
\end{multline*}
%
%
Hence
\begin{equation}   \label{Le4-2}
| E_i(Z_i)|   \leq Z_{i,1}^{'} Z_{i,2}^{'},
\end{equation}
where
\begin{equation}   \label{Le4-22}
     Z_{i,j}^{'}   =\frac{1}{p_i^{r_{i,k_{i,j}} -r_{i,k_{i,j-1}}}}
\Bigg| \sum_{\substack{x_{i,k} \in [0,p_i)\\
k \in [r_{i,k_{i,j-1} }+1,r_{i,k_{i,j}} ]}}
        e\Bigg(  \frac{\hat{m}_{i} \sum_{k=r_{i,k_{i,j-1}}+1}^{r_{i,k_{i,j}} } x_{i,k}p_i^{k-1}}{p_i^{r^{+}_i}}
	\Bigg)  \tilde{\psi}(i ,\hat{\mu}_{i,j} , x_{i,r_{i,j}} ) \Bigg|.
\end{equation}
Let's consider the case $r_{i,k_{i,j}} \leq r_{i,k_{i,j-1}}+1$.
From \eqref{Le2-53} and \eqref{Le3-20}, we get  $|\tilde{\psi}(i,m, x_{i, r_i})| \leq p_i  $,
$ Z_{i,j}^{'} \leq p_i$ and
 $\bar{\mm}_{i,k_{i,j}} \leq p_i^{r_{i,k_{i,j} }-r_{i,k_{i,j-1}}} \leq p_i$.
Taking  $Z_{i,j}^{''}=p_i^2 $ for $ r_{i,k_{i,j}} \leq r_{i,k_{i,j-1}}+1 $,
 we have
\begin{equation} \label{Le2-60}
  Z_{i,j}^{'} \leq  Z_{i,j}^{''}=p_i^2  \leq \frac{p_i^{r_{i,k_{i,j} }-r_{i,k_{i,j-1}}+2}}{\bar{\mm}_{i,k_{i,j}} }
      \leq  \frac{p_i^{3}}{\bar{\mm}_{i,k_{i,j}} }.
\end{equation}
Let's consider the case $r_{i,k_{i,j-1}}+1 < r_{i,k_{i,j}}$.
 Let
\begin{equation}  \nonumber 
     Z_{i,j}^{''}   =\frac{1}{p_i^{r_{i,k_{i,j}} -r_{i,k_{i,j-1}}-2}}
\Bigg| \sum_{\substack{x_{i,k} \in [0,p_i)\\
k \in [r_{i,k_{i,j-1} }+1,r_{i,k_{i,j}} -1 ]}}
        e\Bigg(  \frac{\hat{m}_{i} \sum_{k=r_{i,k_{i,j-1}}+1}^{r_{i,k_{i,j}} -1} x_{i,k}p_i^{k-1}}{p_i^{r^{+}_i}}
	\Bigg)   \Bigg| .
\end{equation}
From \eqref{Le2-53}, \eqref{Le2-01}, \eqref{Le4-2}  and \eqref{Le4-22}, we get  $|\tilde{\psi}(i,m, x_{i, r_i})| \leq p_i  $,  $ Z_{i,j}^{'} \leq Z_{i,j}^{''}$ and
\begin{equation}  \nonumber  
| E_i(Z_i)|   \leq Z_{i,1}^{''} Z_{i,2}^{''}.
\end{equation}
Taking into account \eqref{Le2-20} and \eqref{Le2-60}, we get that in order to prove Lemma 4 it is enough to verify the following inequality
\begin{equation}  \label{Le2-25}
 Z_{i,j}^{''} \leq
 \frac{p_i^4}{\bar{\mm}_{i, k_{i,j}}} \quad  \for \quad r_{i,k_{i,j-1}}+2 \leq r_{i,k_{i,j}}.
\end{equation}
%
%
 Applying  \eqref{Beg-2a},  we derive
\begin{equation} \label{Le2-70}
   Z_{i,j}^{''} \leq p_i \min \Bigg(1,   \frac{1}{2p_i^{r_{i,k_{i,j} }-r_{i,k_{i,j-1}}-1}
		\big\llangle \hat{m}_{i}  /p_i^{r^{+}_i -r_{i,k_{i,j-1} }} \big\rrangle   } \Bigg).
\end{equation}
 If $|\mm_{i, k_{i,j}}| \leq 8$, then we will use the trivial estimate $ Z^{''}_{i,j} \leq p_i \leq p_i^{4}/\bar{\mm}_{i,k_{i,j}}$.
Now let's consider the case $ |\mm_{i, k_{i,j}}| > 8$.
By   \eqref{Le2-70}, we get that in order to prove \eqref{Le2-25} it is enough to verify that
\begin{equation}\label{Le2-80}
   \Big\llangle \frac{\hat{m}_{i} } {p_i^{r^{+}_i -r_{i,k_{i,j-1} }} }
    \Big\rrangle  \geq
     \frac{|\mm_{i,k_{i,j}}| /2}{p_i^{ r_{i,k_{i,j} } -r_{i,k_{i,j-1} }} }.
\end{equation}
From \eqref{Le3-20}, \eqref{Beg-1a} and the previous conditions,  we obtain $r_{i,k_{i,j} } \geq r_{i,k_{i,j-1} }+2 $,
$8 <|\mm_{i, k_{i,j}}| \leq  \frac{1}{2} p_i^{r_{i,k_{i,j} } -r_{i,k_{i,j-1} }}$, $j \in \{1,2\}$ and
\begin{multline*} \nonumber
  \Big\llangle \frac{\hat{m}_{i} } {p_i^{r^{+}_i -r_{i,k_{i,j-1} }} }
    \Big\rrangle = \Big\llangle   \sum_{\ell=1}^{2}   \frac{ \mm_{i, k_{i,\ell}} p_i^{r^{+}_i -r_{i,k_{i,\ell}}} } {p_i^{r^{+}_i -r_{i,k_{i,j-1} }} }
    \Big\rrangle  =
\Big\llangle   \sum_{\ell=1}^{2}   \frac{ \mm_{i, k_{i,\ell}} } {p_i^{r_{i,k_{i,\ell}} -r_{i,k_{i,j-1} }} }
    \Big\rrangle   \\
=\Big\llangle   \sum_{\ell=j}^{2}   \frac{ \mm_{i, k_{i,\ell}} } {p_i^{r_{i,k_{i,\ell}} -r_{i,k_{i,j-1} }} }
    \Big\rrangle   =
\Big\llangle   \frac{ \mm_{i, k_{i,j}} } {p_i^{r_{i,k_{i,j}} -r_{i,k_{i,j-1} }} } +    \sum_{2 \geq \ell>j}  \frac{ \mm_{i, k_{i,\ell}} } {p_i^{r_{i,k_{i,\ell}} -r_{i,k_{i,j-1} }} }
    \Big\rrangle   \\
=\Big\llangle   \frac{ \mm_{i, k_{i,j}} + \varrho_{i,k_{i,j}}} {p_i^{r_{i,k_{i,j}} -r_{i,k_{i,j-1} }} }     \Big\rrangle \quad \with \quad
\varrho_{i,k_{i,j}} =
  \sum_{2 \geq \ell>j}  \frac{ \mm_{i, k_{i,\ell}} } {p_i^{(r_{i,k_{i,\ell}} -r_{i,k_{i,j-1} }) - (r_{i,k_{i,j}} -r_{i,k_{i,j-1} })} }.
\end{multline*}
It is easy to see that $\varrho_{i,k_{i,2}} =0$, and $|\varrho_{i,k_{i,1}}| = |\mm_{i, k_{i,2}} |/p_i^{r_{i,k_{i,2}} -r_{i,k_{i,1} }} \leq 1/2  $.\\
Therefore
\begin{equation}  \nonumber
 \Big\llangle \frac{\hat{m}_{i} } {p_i^{r^{+}_i -r_{i,k_{i,j-1} }} }
    \Big\rrangle = \Big\llangle   \frac{ \mm_{i, k_{i,j}} + \varrho_{i,k_{i,j}}} {p_i^{r_{i,k_{i,j}} -r_{i,k_{i,j-1} }} }     \Big\rrangle
    \geq
    \frac{ |\mm_{i,k_{i,j}}| -2}{p_i^{r_{i,k_{i,j}}-r_{i,k_{i,j-1}}}}
 \geq
     \frac{ |\mm_{i,k_{i,j}}|/2 }{p_i^{r_{i,k_{i,j}}-r_{i,k_{i,j-1}}}}.
\end{equation}
Hence, \eqref{Le2-80}, \eqref{Le2-25} and Lemma 4 are proved.\qed\\  \\
{\bf Lemma 5.}  {\it Let
\begin{equation}  \label{Le4-1}
\DD^{\ast}_{\lambda_1, \lambda_2}= \sum_{ (\br_1,\br_2) \in \UU_{ \lambda_1, \lambda_2 }} \;\;
\sum_{\substack{ \fm_{i,k_{i,j}} \in I_{p_i^{ r_{i,k_{i,j}}  -r_{i,k_{i,j-1}} }}  \\ 0< |m_j|  \leq n^{10}, \; |\fm_{i,j}| \leq n^{10},\; i,j=1,2}} \;
\frac{1}{\bar{m}_1  \bar{m}_{2}}
\;   \prod_{i,j=1}^{2} \frac{1}{\bar{\fm}_{i, j}}  \delta_{p_i^{r^{+}_i}}( \tilde{m}_{i} ) \;,
\end{equation}
%
%
%
\begin{equation}  \nonumber
\DD^{\#}_{\lambda_1, \lambda_2}=
\sum_{\substack{ 0< |m_j|  \leq n^{10}, \; |\fm_{i,j}| \leq n^{10}  \\ i,j=1,2}} \;\;\;   \sum_{ (\br_1,\br_2) \in \UU_{ \lambda_1, \lambda_2 }}\;
\frac{1}{\bar{m}_1  \bar{m}_{2}}
\;   \prod_{i,j=1}^{2}  \frac{1}{\bar{\fm}_{i, j}}  \delta_{p_i^{r^{+}_i}}( \tilde{m}_{i} ),
\end{equation}
where
\begin{equation}  \label{Le3-18}
 \tilde{m}_{i}   = \tilde{m}_{i}( m_{1}, m_{2}, \fm_{i,1}, \fm_{i,2}  )=
 \sum_{j=1}^{2}  (\fm_{i,j} + m_{j} M_{i, \br_{j}})p_i^{ r^{+}_i -  r_{i,j}  } .
\end{equation}
Then}
\begin{equation}  \nonumber
  \Big| \int_{[0,1]^2}  \bD_{\lambda_1, \lambda_2,N}   d \bx \Big|   \ll   1 +   \DD^{\ast}_{\lambda_1, \lambda_2}  \quad \ad \quad
       \quad     \DD^{\ast}_{\lambda_1, \lambda_2}  \leq \DD^{\#}_{\lambda_1, \lambda_2}   .
\end{equation} \\
%
{\bf Proof.}
From \eqref{Le3-9a} and \eqref{Le3-20},  we derive
\begin{multline} \label{Le3-21}    
  \tilde{m}_{i}( m_{1}, m_{2}, \mm_{i,1}, \mm_{i,2}  )=
 \sum_{j=1}^{2}  (\mm_{i,j} + m_{j} M_{i, \br_{j}})p_i^{ r^{+}_i -  r_{i,j}} \\
 \equiv  \hat{m}_i +\sum_{j=1}^{2}   m_{j} M_{i, \br_{j}} p_i^{ r^{+}_i -  r_{i,j}} \equiv
 0    \mod p_i^{r^{+}_i}.
\end{multline}
By \eqref{Le3-20}, $ \mm_{i,k_{i,j}} \in I_{p_i^{ r_{i,k_{i,j}}  -r_{i,k_{i,j-1}} }}$. Together with \eqref{Le4-0}, this implies that
\begin{multline} \label{Le3-ab}    
 \Big| \int_{[0,1]^2}  \bD_{\lambda_1, \lambda_2,N}   d \bx \Big|   \ll
\sum_{ (\br_1,\br_2) \in \UU_{ \lambda_1, \lambda_2 }} \;\;\;
\sum_{  m_j \in I^{*}_{P_{\br_j}}, \;  j=1,2} \;
\frac{ 1 }{\bar{m}_1  \bar{m}_{2}}   \\
 \times \sum_{ \fm_{i,k_{i,j}} \in I_{p_i^{ r_{i,k_{i,j}}  -r_{i,k_{i,j-1}} }}, \;  i,j=1,2} \;
\;   \prod_{i,j=1}^{2} \frac{1}{\bar{\fm}_{i, j}}  \delta_{p_i^{r^{+}_i}}( \tilde{m}_{i} )
 .
\end{multline}

 Bearing in mind \eqref{Le3-18}, we get that the part of the right hand of \eqref{Le3-ab},  satisfying  the condition $|\fm_{i,j}| >n^{10}$ for some
$i \in \{ 1,2 \},\; j \in \{ 1,2 \} $, is equal to $O(1)$.
Therefore
\begin{multline*}  
 \Big| \int_{[0,1]^2}  \bD_{\lambda_1, \lambda_2,N}   d \bx \Big|   \ll 1 +	
   \sum_{ (\br_1,\br_2) \in \UU_{ \lambda_1, \lambda_2 }}
\sum_{\substack{ m_j \in I^{*}_{P_{\br_j}} \\ j=1,2}}
\sum_{ \fm_{i,k_{i,j}} \in I_{p_i^{ r_{i,k_{i,j}}  -r_{i,k_{i,j-1}} }}, \;  i,j=1,2} \; \frac{1 }{\bar{m}_1 \bar{m}_{2}}   \\
   \times  \prod_{i,j=1}^{2}  \frac{\Delta(|\fm_{i,j}| \leq n^{10})}{\bar{\fm}_{i,j}}
   \delta_{p_i^{r^{+}_i}} ( \tilde{m}_{i}  )
			.
\end{multline*}
%
Let
\begin{equation*} 
\theta_1(\bm) =\Delta(\max_j |m_j| \leq n^{10}), \qquad \theta_2(\bm) =1-\theta_1(\bm).
\end{equation*}
We have
\begin{multline}  \label{L5a-6}
 \Big| \int_{[0,1]^2}  \bD_{\lambda_1, \lambda_2,N}   d \bx \Big|  \ll 1 +	\check{\DD}_{\lambda_1, \lambda_2,1}+ \check{\DD}_{\lambda_1, \lambda_2,2}, \;\;  \with \\
   \check{\DD}_{\lambda_1, \lambda_2,\nu}=
\sum_{ (\br_1,\br_2) \in \UU_{ \lambda_1, \lambda_2 }} \;
     \sum_{\substack{ \fm_{i,k_{i,j}} \in I_{p_i^{ r_{i,k_{i,j}}  -r_{i,k_{i,j-1}} }}    \\ |\fm_{i,j}|  \leq n^{10}, i,j=1,2}}\;\; \sum_{m_j \in I^{*}_{P_{\br_j}}, \;   j=1,2}
   \prod_{j=1}^{2} \frac{1 }{\bar{m}_j} \\
  \times
      \prod_{i =1 }^2 \frac{ \theta_{\nu}(\bm)}{\bar{\fm}_{i,j}}  \delta_{p_i^{r^{+}_i}} ( \tilde{m}_{i}  ) , \quad \nu=1,2.
\end{multline}	
From \eqref{Le4-1}, we obtain  $ \DD^{\ast}_{\lambda_1, \lambda_2} \leq \DD^{\#}_{\lambda_1, \lambda_2} $ and  $ \check{\DD}_{\lambda_1, \lambda_2,1} \leq \DD^{\ast}_{\lambda_1, \lambda_2} $.
Hence, in order to get the assertion of Lemma 4, it is enough to prove that
\begin{equation} \label{L5a-7}
     \check{\DD}_{\lambda_1, \lambda_2,2} \ll 1.
\end{equation}
We consider $\check{\DD}_{\lambda_1, \lambda_2,2}$  (the case of $\theta_{2}(\bm) =1) $.
Let $|m_{j_0}| \geq n^{10} $  with some $j_0 \in [1,2]$. Let $j_0=1$. The case $j_0=2$ is similar.
 By  \eqref{Beg-26}, we get that $\max ( r_{1,1},r_{2,1}) > V=p_1^4 p_2^4 [\log_2^4 n]$. Therefore $ r_{i_0,1} \geq V$  with some $i_0 \in [1,2]$.
In view of   \eqref{L5a-6}, we have  that \eqref{L5a-7} ensue from the following inequality
\begin{equation} \label{Le5-7}
W=W(\br_1,\br_2, m_2,\fm_{1,1},\fm_{1,2},\fm_{2,1},\fm_{2,2}  ):= \sum_{p_1^n p_2^n \geq |m_{1}| \geq n^{10} } \frac{1}{\bar{m}_{1}}  \delta_{p_{i_0}^{r^{+}_{i_0}}}
  ( \tilde{m}_{{i_0}}  ) \ll  n^{-9}.
\end{equation}
Now we will prove \eqref{Le5-7}. We fix $\fm_{i,j}\; (i, j \in [1, 2])$, $\br_1,\br_{2}$ and $ m_2$.
Let
\begin{equation} \nonumber
K_1=K_1(\br_1,\br_2, m_2,\fm_{1,1},\fm_{1,2},\fm_{2,1},\fm_{2,2}  ):=- \sum_{j=1}^{2}    \fm_{{i_0},j}  p_{i_0}^{ r^{+}_{i_0} -   r_{i_0,j}}
- m_{2} M_{i_0, \br_{2}}  p_{i_0}^{ r^{+}_{i_0} -   r_{i_0,2}}
\end{equation}
and let
\begin{multline} \label{Le5-77}
W=W_1+W_2,   \quad \with  \quad W_{\nu}= \sum_{p_1^n p_2^n \geq |m_{1}| \geq n^{10} } \frac{1}{\bar{m}_{1}}  \delta_{p_{i_0}^{r^{+}_{i_0}}}
  ( \tilde{m}_{{i_0}}  ) \tilde{\theta}_{\nu}(K_1), \quad \where   \\
  \tilde{\theta}_{1}(K_1)=\delta_{p_{i_0}^{r^{+}_{i_0} -r_{i_0,1} }}
  ( K_1 ),
    \;\;\;\; \tilde{\theta}_{2}(K_1) =1-\tilde{\theta}_{1}(K_1).
\end{multline}
By   \eqref{Le3-18} and \eqref{Le5-7}, we get
\begin{equation} \nonumber
   \tilde{m}_{{i_0}} = \sum_{j=1}^{2}   ( \fm_{{i_0},j} + m_{j} M_{i_0, \br_{j}} ) p_{i_0}^{ r^{+}_{i_0} -   r_{i_0,j}}
\equiv 0 \mod p_{i_0}^{ r^{+}_{i_0}}.
\end{equation}
Hence
\begin{equation} \label{Le5-78}
 K_1 \equiv m_{1}  M_{i_0, \br_{1}}  p_{i_0}^{ r^{+}_{i_0} -   r_{i_0,1}}  \mod p_{i_0}^{ r^{+}_{i_0}},
 \quad \tilde{\theta}_{1}(K_1) =1 \quad \ad \quad W_2=0.
\end{equation}
Let's consider $W_1$.  
 Let $K_2 = K_1    p_{i_0}^{ -r^{+}_{i_0} +   r_{i_0,1}}$.
%
Then  $ m_{1}  M_{i_0, \br_{1}}  \equiv K_2 \mod p_{i_0}^{ r_{i_0,1}} $.
According to \eqref{Beg-3}, $ M_{i_0, \br_{1}}
  p_{i_0^{'} }^{r_{i_0^{'},1} }  \equiv 1 \mod p_{i_0}^{r_{i_0,1}}$, $ i_0^{'} \in \{1,2\} \setminus \{i_0\}$.\\
Let $ K_3= K_2 p_{i_0^{'} }^{r_{i_0^{'},1} } $.
Therefore  $   m_{1} \equiv K_3    \mod p_{i_0}^{r_{{i_0},1}}  $.
Taking into account that $ r_{i_0,1} \geq V \geq \log^4_2 n , \;\; p_{i_0}^{r_{i_0,1}} \gg n^{20}  $, from \eqref{Le5-7}, we have
\begin{equation*}  
W \ll  \sum_{\ell \in \ZZ} \frac{1}{|K_3 +\ell p_{i_0}^{r_{i_0,1}}|}
 \Delta \Big( p_1^n p_2^{ n}  \geq |K_3 + \ell p_{i_0}^{r_{i_0,1}}|  \geq n^{10} \Big)
 \ll \sum_{\ell =1}^{ p_1^n p_2^{ n}} \frac{1}{\ell  n^{10}}
   \ll  n^{-9}.
\end{equation*}
Bearing in mind \eqref{Le5-77} and \eqref{Le5-78}, we have that \eqref{Le5-7}, \eqref{L5a-7} and  Lemma 5 are proved. \qed  \\

By   \eqref{Beg-28}, \eqref{Beg-90} and Lemma 3, we get that in order to prove the Theorem it is enough to verify that
\begin{equation}   \label{Le4-30}
  \DD^{\ast}_{\lambda_1, \lambda_2}     \ll  n^2, \quad \for \quad \lambda_i \in \{ 1,2 \}, \; i=1,2.
\end{equation}  \\

 {\bf 3. \;\;  The main lemmas.}\\
{\bf 3.1. \;\; Linear forms in $p$-adic logarithm.}
Let $\tilde{\alpha}_1,...,\tilde{\alpha}_{\tilde{n}}$  be non-zero algebraic numbers and $K$ be a number
field containing $\tilde{\alpha}_1,...,\tilde{\alpha}_{\tilde{n}}$ with $d=[K: \QQ]$. Denote by $\fd$  a prime ideal
of the ring $\fO_K$ of integers in $K$, lying above the prime number $p$, and by $\ff_{\fd}$  the residue class degree of $\fd$. For
$\gamma \in K$, $\gamma \neq 0$, write ${\rm ord}_{\fd}(\gamma)$ for the exponent to which $\fd$ divides the principal fractional ideal generated by $\gamma $ in $ K$.
Define
\begin{equation*}  
    h^{'}(\tilde{\alpha}_j) = \max( h_0(\tilde{\alpha}_j), \ff_{\fd} (\log p)/d)
		\qquad  (1 \leq j \leq \tilde{n}),
\end{equation*}
where $h_0(\gamma)$ denotes the absolute logarithmic Weil height of an algebraic
number $\gamma$, i.e.,
\begin{equation*}  
    h_0(\gamma) = \tilde{k}^{-1} \Big(  \log a_0  + \sum_{i=1}^{\tilde{k}} \log \max(1, |\gamma^{(i)}|) \Big) ,
\end{equation*}
where the minimal polynomial for $\gamma$ is
\begin{equation*}  
  a_0x^{\tilde{k}} +   a_1 x^{{\tilde{k}}-1} + \cdots + a_{\tilde{k}} =a_0(x-\gamma^{(1)}) \cdots (x-\gamma^{({\tilde{k}})}), \qquad a_0 >0.
\end{equation*}

{\bf Theorem A.} ([Yu, Theorem 1], [Bu, Theorem 2.9]) {\it Let $\tilde{\Xi} = \tilde{\alpha}_1^{b_1} \cdots \tilde{\alpha}_{\tilde{n}}^{b_{\tilde{n}}} -1\\ \neq 0 $. Suppose that
\begin{equation*}  
    {\rm ord}_{\fd}(\tilde{\alpha}_j) =0 \qquad (1 \leq j \leq \tilde{n}).
\end{equation*}
Then there exists a constant $C_0$, depending only on $\tilde{n}, d$ and $\fd$, such that
\begin{equation*}  
    {\rm ord}_{\fd}(\tilde{\Xi})  < C_0 h^{'}(\tilde{\alpha}_1) \cdots h^{'}(\tilde{\alpha}_{\tilde{n}}) \log \tilde{B},
\end{equation*}
where}
\begin{equation*}  
\tilde{B} = \max(|b_1|,...,|b_{\tilde{n}}| ,3), \quad b_i \in \ZZ.
\end{equation*}

We will use Theorem A with $\tilde{n}=2,\;{\tilde{k}} = d=1,\; \fd=p_i$,
$ \tilde{\alpha}_1 \in \{p_1,p_2\}$, $\tilde{\alpha}_2 =l_1/l_2$, $l_1,l_2, b_1,b_2 \in \ZZ$.\\

{\bf Corollary.}  {\it Let $0<|l_i| \leq n^{{10}}$,  $  {\rm ord}_{p_i}(l_1/l_2)=0$,  $i^{'} \in \{ 1,2 \} \setminus \{ i \}$,
      $(i =1,2)$ and $\tilde{B} =n$. Then there exists a constant $C_1>0$ such that }
\begin{equation*}  
   {\rm ord}_{p_i}\Big( \Xi_i \Big)  < C_1 \log_2 n \log_2 \tilde{B} =C_1 \log_2^2n, \;\;\; \with \;\;\;
   \Xi_i =(l_1 / l_2) p_{i^{'}}^{b_{i^{'} }} -1 \neq 0, \;\;i=1,2.
\end{equation*}\\

 {\bf 3.2. \;\;  The applications of Theorem A.}\\ \\
%
%
%
We denote by $ \DD^{\ast,j}_{\lambda_1, \lambda_2} $ a part of $ \DD^{\ast}_{\lambda_1, \lambda_2} $
(respectively $ \DD^{\#,j}_{\lambda_1, \lambda_2} $ a part of $ \DD^{\#}_{\lambda_1, \lambda_2} $)
(see \eqref{Le4-1})  satisfying the condition number $j \in \{1,...,5 \}$:
\begin{multline}  \label{Le7-0}
 \DD^{\ast,j}_{\lambda_1, \lambda_2} = \sum_{ (\br_1,\br_2) \in \UU_{ \lambda_1, \lambda_2 }} \;\;
\sum_{\substack{ \fm_{i,k_{i,j}} \in I_{p_i^{ r_{i,k_{i,j}}  -r_{i,k_{i,j-1}} }}  \\ 0< |m_j|  \leq n^{10}, \; |\fm_{i,j}| \leq n^{10}}} \;
\frac{1}{\bar{m}_1  \bar{m}_{2}}
\;   \prod_{i,j=1}^{2} \frac{1}{\bar{\fm}_{i, j}}  \delta_{p_i^{r^{+}_i}}( \tilde{m}_{i} )   \Delta(case \; j ),\\
 \DD^{\#,j}_{\lambda_1, \lambda_2} = \sum_{ (\br_1,\br_2) \in \UU_{ \lambda_1, \lambda_2 }} \;\;
\sum_{ 0< |m_j|  \leq n^{10}, \; |\fm_{i,j}| \leq n^{10}} \;
\frac{1}{\bar{m}_1  \bar{m}_{2}}
\;   \prod_{i,j=1}^{2} \frac{1}{\bar{\fm}_{i, j}}  \delta_{p_i^{r^{+}_i}}( \tilde{m}_{i} )   \Delta(case \; j ).
\end{multline}
 We will prove that $ \DD^{\ast,j}_{\lambda_1, \lambda_2} \ll n^2 $ (see  \eqref{Le4-30}).
 We consider the case $j \in [1,3]$ in Lemma 6 and the case  $j \in [4,5]$ in Lemma 7.
 Proofs of Lemma~6 and Lemma 7 are similar. To make the article more readable, we do not combine Lemma~6 and Lemma 7. \\  \\
{\bf Lemma 6.} {\it  The estimate \eqref{Le4-30} is true for  $ \lambda_{1} = \lambda_{2}=0$.}\\ \\
{\bf Proof.} We will consider \eqref{Le4-30} for $n \geq n_0$, with $ n_0 = [2^ {4C_1+20}]p_1^4 p^4_2$, where $C_1$ is defined in the Corollary.
From     \eqref{Beg-26} and  \eqref{Beg-40}, we have for $n \geq n_0$, $ i,j \in \{1,2\}$, that
\begin{multline}   \label{Le5-4}
\hat{r}_j=  \max(r_{1,j}, r_{2,j})> V, \quad \quad r^{+}_i -  r^{-}_i > \VV \;\;\for \;\;\lambda_{i}=1, \quad \quad   r^{+}_i -  r^{-}_i \leq \VV \;\;\for \; \; \lambda_{i}=0,  \\
  r^{+}_i =\max(  r_{i,1}, r_{i,2} ), \quad r^{-}_i =\min(  r_{i,1}, r_{i,2} ), \quad V =p_1^4 p_2^4 [\log_2^4 n]  \geq 4(C_1+10) \log_2^2 n, \\
   V/(4 \log_2(p_1p_2))  > \VV = [\log_2^3 n]\geq 4C_1 \log_2^2 n.
\end{multline}
  {\it Case} 1: $ \lambda_{1} = \lambda_{2}=0$ and $r^{+}_{i_0} \leq V$ for some  $i_0 \in \{1,2\}$.  Using  \eqref{Le7-0}, we derive

\begin{multline}  \label{Le4-98}
 \DD^{\#,1}_{0, 0} = \sum_{ (\br_1,\br_2) \in \UU_{ 0, 0 }} \;\;
\sum_{ 0< |m_j|  \leq n^{10}, \; |\fm_{i,j}| \leq n^{10},\; i,j=1,2} \;
\frac{1}{\bar{m}_1  \bar{m}_{2}}
\;   \prod_{i,j=1}^{2} \frac{1}{\bar{\fm}_{i, j}}  \delta_{p_i^{r^{+}_i}}( \tilde{m}_{i} )   \Delta(r^{+}_{i_0} \leq V )\\
 \ll
\sum_{ (\br_1,\br_2) \in \UU_{ 0, 0 }}\; \log^6 n \;
  \Delta(r^{+}_{i_0} \leq V )  \leq
\sum_{\substack{ r^{+}_i -  r^{-}_i \leq \VV, r^{+}_{i_0} \leq V\\
 r_{i,j} \in [1,n], \; i,j=1,2}} \log^6 n
   \ll nV \VV^2 \log^6 n    \ll n^2.
\end{multline}
%
In the following, we consider the case $\min (r^{+}_1,r^{+}_2) > V $.
Let
\begin{equation}   \label{Le5-8}
    g_{i}  :=  \sum_{j=1}^2 \fm_{i,j} p_{i}^{r^{+}_{i} - r_{i,j} }, \quad
    i^{'} \equiv  i+1 \mod 2, \;\;   i^{'} \in \{ 1,2 \} , \;\; i=1,2.
\end{equation}
%
%
By \eqref{Beg-3}, \eqref{Le3-18} and \eqref{Le7-0}, we obtain that
\begin{multline}   \label{Le5-2}            
   \tilde{m}_{i}  =\sum_{j=1}^{2}  (\fm_{i,j} + m_{j} M_{i, \br_{j}})p_i^{ r^{+}_i -  r_{i,j}  }
    \mod p_i^{r^{+}_i}, \quad  M_{i,\br_j}  p_{i^{'}}^{r_{i^{'},j} } \equiv 1 \mod p_i^{r_{i,j}},\\
   p_{i^{'}}^{r^{+}_{i^{'}}} \tilde{m}_{i} =  g_{i} p_{i^{'}}^{r^{+}_{i^{'}}}  + \sum_{j=1}^2 m_{j}  M_{i,\br_j}   p_{i^{'}}^{r^{+}_{i^{'}}}  p_{i}^{r^{+}_{i} - r_{i,j} }
     , \quad
       M_{i,\br_j}  p_{i^{'}}^{ r^{+}_{i^{'}} } \equiv   p_{i^{'}}^{ r^{+}_{i^{'}}- r_{i^{'},j} }  \mod p_i^{r_{i,j}},\\
  M_{i,\br_j} p_{i^{'}}^{ r^{+}_{i^{'}} }    p_i^{r^{+}_i  -r_{i,j}}  \equiv
  p_{i^{'}}^{ r^{+}_{i^{'}} - r_{i^{'},j} } p_i^{r^{+}_i  -r_{i,j}} \mod p_i^{r^{+}_i}.
\end{multline}
%
Hence
\begin{equation}    \label{Le5-10}
    g_{i} p_{i^{'}}^{r^{+}_{i^{'}}} + \sum_{j=1}^2 m_{j}  p_{i}^{r^{+}_{i} - r_{i,j} } p_{i^{'}}^{r^{+}_{i^{'}} - r_{i^{'},j} } \equiv 0 \mod p_i^{r^{+}_i}, \quad i=1,2.
\end{equation}
%
%
{\it Case} 2: $ \lambda_{1} = \lambda_{2}=0$, $\quad g_1=g_2=0$ and $\min (r^{+}_1,r^{+}_2) > V $.
From \eqref{Le5-10}, we have
\begin{equation}  \label{Le5-14}
    \sum_{j=1}^2 m_{j}  p_{i}^{r^{+}_{i} - r_{i,j} } p_{i^{'}}^{r^{+}_{i^{'}} - r_{i^{'},j} } \equiv 0  \mod p_i^{r^{+}_i} .
\end{equation}
Bearing in mind  \eqref{Le5-4}, and that $0 < |m_j|  \leq  n^{10}$ ($j=1,2)$,
$ \lambda_{1} = \lambda_{2}=0$,
$ V/(4\log_2(p_1p_2))  $  $\geq  \VV = [\log_2^3 n] \geq r^{+}_i -  r^{-}_i$,
$r^{+}_i >V$ $(i=1,2)$, we get
\begin{equation} \label{Le5-30}
     \log_2 \big(\big| \sum_{j=1}^2 m_{j}  p_{i}^{r^{+}_{i} - r_{i,j} }  p_{i^{'}}^{r^{+}_{i^{'}} - r_{i^{'},j} } \big|\big) \leq 20 \log_2 n + \VV \log_2 (p_1 p_2) \leq  V -\VV <r^{+}_i -\VV.
\end{equation}
%
Hence the congruence   \eqref{Le5-14} is equality  :
\begin{equation}  \nonumber 
     \sum_{j=1}^2 m_{j}  p_{i}^{r^{+}_{i} - r_{i,j} } p_{i^{'}}^{r^{+}_{i^{'}} - r_{i^{'},j} } =0, \quad i=1,2.
\end{equation}
 Let $\mu_3 \geq 1$  be the greatest common divisor of $ m_1$ and $m_2$, and let $ m_1^{*}=m_1/\mu_3$, $ m_2^{*}=m_2/\mu_3$,   $ ( m_1^{*}, m_2^{*})=1 $. Then
\begin{equation}  \nonumber 
 \sum_{j=1}^2 m_j^{*} p_{i}^{r^{+}_{i} - r_{i,j} } p_{i^{'}}^{r^{+}_{i^{'}} - r_{i^{'},j} } =0.
\end{equation}
It is easy to verify that
 $ |m_1^{*}  m_2^{*}| =   p_1^{r^{+}_{1}- r^{-}_1} p_2^{r^{+}_{2} - r^{-}_2} $ and $ |m_j^{*}|  \in \{1, p_1^{r^{+}_{1}- r^{-}_1},
   p_2^{r^{+}_{2} - r^{-}_2}, \\  p_1^{r^{+}_{1}- r^{-}_1} p_2^{r^{+}_{2} - r^{-}_2} \}$, $j=1,2$.

 Let $\mu_i \geq 1$  be the greatest common divisor of $ \fm_{i,1}$ and $\fm_{i,2}$.
Taking into account that $ g_1=g_2=0$,    we have from \eqref{Le5-8} that
 $ \fm_{i,1}  = \fm_{i,2} =0 $ or
  $ |\fm_{i,1}\fm_{i,2}| =\mu_i^2  p_i^{r^{+}_{i}- r^{-}_i} $ and  $ |\fm_{i,j}| /\mu_i
  \in\{1,  p_i^{r^{+}_{i}- r^{-}_i}\} $ , $i,j=1,2$.

Substituting the above relations into \eqref{Le7-0}, we get
\begin{equation}  \nonumber
\DD^{\#,2}_{0, 0}  \leq
\sum_{\substack{  0< |m_j|  \leq n^{10}, \; |\fm_{i,j}| \leq n^{10} \\ i,j=1,2}} \;
  \sum_{ \substack{ (\br_1,\br_2) \in \UU_{ 0, 0 }   \\ \min (r^{+}_1,r^{+}_2) > V   }}\;
\frac{1}{\bar{m}_1  \bar{m}_{2}}
\;   \prod_{i,j=1}^{2} \frac{\Delta(g_i=0)}{\bar{\fm}_{i, j}} \delta_{p_i^{r^{+}_i}}( \tilde{m}_{i} )
\end{equation}
\begin{equation}  \label{Le4-90}
\ll
\sum_{\substack{  \mu_i \geq 1   \\ i=1,2,3}} \;
  \sum_{\substack{  r_{i,j} \in [1,n]   \\ i,j=1,2}} \;
\frac{1}{\mu_1^2 \mu_2^2  \mu_3^2    }
\;    p_1^{-r^{+}_{1}+ r^{-}_1} p_2^{-r^{+}_{2}+ r^{-}_2} \ll    \sum_{ \substack{ r_{i,j} \in [1,n]   \\ i,j=1,2}}
 p_1^{-r^{+}_{1}+ r^{-}_1} p_2^{-r^{+}_{2}+ r^{-}_2}
 \ll n^2.
\end{equation}\\
{\it Case} 3: $ \lambda_{1} = \lambda_{2}=0$, $\quad |g_1| + |g_2| >0$ and $\min (r^{+}_1,r^{+}_2) > V $.
   Let  $g_1  \neq 0$. The proof for the case  $g_2  \neq 0$ is similar.
By  \eqref{Le5-10}, we get
\begin{equation} \label{Le6-50}
    g_{1} p_2^{r^{+}_{2}} + \xi \equiv 0 \mod p_1^{r^{+}_{1}},
     \quad \with \quad \xi=\sum_{j=1}^2 m_{j} p_{1}^{r^{+}_{1} - r_{1,j} } p_{2}^{r^{+}_{2} - r_{2,j} }, \quad r^{+}_1 > V.
\end{equation}
Bearing in mind that $ \lambda_{1} = \lambda_{2}=0$, we get
from \eqref{Le5-4}
that $ \rho_i := r^{+}_{i} -r^{-}_{i}
\leq \VV$ for $i=1,2$.
We fix $m_{1}, m_{2}, \fm_{i,j}$ and $ \rho_i $, with  $i,j=1,2$.
Let $\beta_1 =\ord_{p_1}(g_{1})$ and let $g_{1}=g_{1}^{'} p_{1}^{\beta_1}$, $ (p_{1},g_{1}^{'} )  =1$.
Similarly to \eqref{Le5-30}, we obtain from \eqref{Le5-8} and \eqref{Le6-50}  that
\begin{equation} \nonumber  
  \beta_1 \leq
  \log_2|\sum_{j=1}^2 \fm_{1,j} p_{1}^{r^{+}_{1} - r_{1,j} }|)
  \leq 20 \log_2 n + \VV \log_2 (p_1 p_2) \leq  V -\VV < r^{+}_{1}-\VV  .
\end{equation}
Let $\beta_2 =\ord_{p_1} (\xi)$.
From  \eqref{Le6-50},  we have that $\beta_2 =\beta_1$    and
\begin{equation}  \label{Le6-70}
    p_2^{r^{+}_{2}} \equiv - \xi^{'} /g_1^{'} \mod p_{1}^{\VV},  \quad \with \quad \xi^{'} = \xi /p_1^{\beta_2}.
\end{equation}
Suppose that \eqref{Le6-70} have two solutions  $r^{+'}_{2}$ and $r^{+''}_{2}$. Then
\begin{equation} \nonumber 
   p_2^{z}  \equiv 1 \mod p_1^{\VV},  \quad \with \quad  z= r^{+'}_{2} -r^{+''}_{2}
   \in [-n,n].
\end{equation}
Hence $\ord_{p_1}(  p_2^{z} - 1 )  \geq \VV = [\log_2^3 n] > C_1 \log_2^2 n$.
By the Corollary, we get that the number of solutions of the above congruence is only one $(z=0)$.
Taking into account that $\lambda_1 =\lambda_2 =0$, we get from \eqref{Le5-4}
that $  0 \leq r^{+}_{i} -r^{-}_{i} \leq \VV$, $i=1,2$.
Therefore, the number of  $(\br_1, \br_2)$ satisfying  \eqref{Le6-70} is less than
\begin{equation}  \nonumber
\#\{ (r^{+}_{1},\rho_1,\rho_2)\; |\; 1 \leq r^{+}_{1} \leq n, 0 \leq \rho_1,\rho_2 \leq \VV  \}=n(\VV+1)^2.
\end{equation}
Using  \eqref{Le7-0}, we derive
\begin{equation}  \label{Le4-94}
 \DD^{\#,3}_{0, 0}\ll
\sum_{\substack{ 0< |m_j|  \leq n^{10}, \; |\fm_{i,j}| \leq n^{10} \\ i,j=1,2}} \;
\frac{n\VV^2}{\bar{m}_1  \bar{m}_{2}}
\;   \prod_{i,j=1}^{2} \frac{1}{\bar{\fm}_{i, j}}  \ll n \VV^2 \log^6 n \ll n^2.
\end{equation}
%
%
From \eqref{Le4-1}, \eqref{Le7-0}, \eqref{Le4-98}, \eqref{Le4-90} and  \eqref{Le4-94}, we obtain $\DD^{\ast}_{0, 0} \leq \DD^{\#}_{0, 0}
  \leq \DD^{\#,1}_{0,0} + \DD^{\#,2}_{0,0} + \DD^{\#,3}_{0,0} \ll n^2$.
Hence Lemma 6 is proved. \qed \\ \\
{\bf Lemma 7.}  {\it  The estimate \eqref{Le4-30} is true for  $\max_i \; \lambda_{i}=1$.}  \\ \\
{\bf Proof.} 
 Let's consider the case  $ \lambda_{1}=1$. The proof for the case  $ \lambda_{2}=1$ is similar.
 We will consider $ \DD^{\ast,j}_{\lambda_1, \lambda_2}$.
  By  \eqref{Beg-26}, \eqref{Beg-40}, \eqref{Le7-0} and \eqref{Le5-4}, we get that
 \begin{multline}\label{Le6-1}
  \hat{r}_j=\max(r_{1,j}, r_{2,j})> V, \quad
               r_{i,k_{i,2}} = r^{+}_i =\max_j  r_{i,j}, \;\; r_{i,k_{i,1}} = r^{-}_i =\min_j  r_{i,j},\\
                 r^{+}_1 -  r^{-}_1  > \VV = [\log_2^3 n], \quad \ad \quad
        \fm_{i,k_{i,j}} \in I_{p_i^{ r_{i,k_{i,j}}  -r_{i,k_{i,j-1}} }}, \quad i,j =1,2.
\end{multline}
Using  \eqref{Le5-2}  and \eqref{Le6-1}, we obtain for $i \in \{1,2  \}$ that
\begin{equation} \label{Le6-2}    
   \tilde{m}_{i}   \equiv   \fm_{i,k_{i,2}} + m_{k_{i,2}} M_{i, \br_{k_{i,2}}}+ (\fm_{i,k_{i,1}} + m_{k_{i,1}} M_{i, \br_{k_{i,1}}})
 p_i^{ r^{+}_i -  r_{i,k_{i,1}}  } \equiv
 0    \mod p_i^{r^{+}_i},
\end{equation}
\begin{equation} \nonumber 
\fm_{i,k_{i,2}} + m_{k_{i,2}} M_{i, \br_{k_{i,2}}}
 \equiv
 0    \mod  p_i^{ r^{+}_i -  r^{-}_i  }, \quad  M_{i, \br_{k_{i,2}}}  p_{{i^{'}}}^{r_{{i^{'},k_{{i,2}} }}} \equiv  1 \mod p_i^{ r^{+}_i  }, \;  r^{+}_i= r_{i,k_{i,2}}.
\end{equation}
Hence
\begin{equation} \label{Le6-3}    
p_{{i^{'}}}^{r_{{i^{'},k_{{i,2}}} }}   \fm_{i,k_{i,2}} + m_{k_{i,2}} +p_{{i^{'}}}^{r_{{i^{'},k_{{i,2}}} }}(\fm_{i,k_{i,1}} + m_{k_{i,1}} M_{i, \br_{k_{i,1}}})
 p_i^{ r^{+}_i -  r_{i,k_{i,1}}  } \equiv
 0    \mod p_i^{r^{+}_i}.
\end{equation}
From \eqref{Le6-1} and \eqref{Le6-2}, we have  $r^{+}_1 -  r^{-}_1  >  \VV$ and
\begin{equation} \label{Le6-4}    
 p_2^{r_{2,k_{1,2}}}    \fm_{1,k_{1,2}} + m_{k_{1,2}}
 \equiv
 0    \mod  p_1^{\VV }.
\end{equation}
Let $\beta_3 =\ord_{p_1}(m_{k_{1,2}})$ and let $m_{k_{1,2}}^{*} =m_{k_{1,2}} p_{1}^{-\beta_3}$, $ (p_{1},m_{k_{1,2}}^{*} )  =1$. Bearing in mind that $ 0 <\max(|m_j |, |\fm_{i,j}|) \leq n^{10}$,
we get $\beta_3 \leq  \VV/2$.    Let $\fm_{1,k_{1,2}}^{*}=\fm_{1,k_{1,2}} p_{1}^{-\beta_3}$.
 Then $     p_2^{r_{2,k_{1,2}}}  \fm_{1,k_{1,2}}^{*} / m_{k_{1,2}}^{*} + 1
 \equiv
 0    \mod  p_1^{ [\VV/2 ]  } $ and
\begin{equation} \nonumber 
 \ord_{p_1}(     p_2^{r_{2,k_{1,2}}}  \fm_{1,k_{1,2}}^{*} / m_{k_{1,2}}^{*} + 1 )  \geq [\VV
/2]  = [ [\log_2^3 n]/2] \geq C_1 \log_2^2 n , \quad n \geq n_0.
\end{equation}
Applying the Corollary, we get that the above congruence is equality
\begin{equation} \label{Le6-10}    
    p_2^{r_{2,k_{1,2}}}  \fm_{1,k_{1,2}}^{*}  =- m_{k_{1,2}}^{*} \quad \ad \quad
    p_2^{r_{2,k_{1,2}}}  \fm_{1,k_{1,2}}  =- m_{k_{1,2}}, \quad \beta_3 =\ord_{p_1}(m_{k_{1,2}}).
\end{equation}
Taking into account that $\max(|m_j |, |\fm_{i,j}|) \leq n^{10}$, we get

\begin{equation} \label{Le6-11}    
 r_{2,k_{1,2}} \leq
  \log_{p_2} (\max(|m_1|, |m_2|)) \leq 10 \log_2 n  \leq \VV .
\end{equation}
%
%
{\it Case} 4 : $\quad  \lambda_{1} = \lambda_{2} =  1$.
  Repeating \eqref{Le6-1} - \eqref{Le6-10} for $\lambda_2 =1$, we get
\begin{equation} \label{Le6-12}    
    p_1^{r_{1,k_{2,2}}}  \fm_{2,k_{2,2}}^{*}  =- m_{k_{2,2}}^{*}  \quad \ad \quad
    p_1^{r_{1,k_{2,2}}}  \fm_{2,k_{2,2}}  =- m_{k_{2,2}},
\end{equation}
with  $\beta_4 =\ord_{p_2}(m_{k_{2,2}})$, $m_{k_{2,2}}^{*}=m_{k_{2,2}} p_{2}^{-\beta_4}$, $ (p_{2},m_{k_{2,2}}^{*} )  =1$,  $\fm_{2,k_{2,2}}^{*} = \fm_{2,k_{2,2}} p_{2}^{-\beta_4}$.\\
Now substituting  \eqref{Le6-10} and \eqref{Le6-12} into \eqref{Le6-3}, we obtain
\begin{equation} \label{Le6-13}    
p_{{i^{'}}}^{r_{{i^{'},k_{{i,2}}} }}(\fm_{i,k_{i,1}} + m_{k_{i,1}} M_{i, \br_{k_{i,1}}})
 p_i^{ r^{+}_i -  r_{i,k_{i,1}}  } \equiv
 0    \mod p_i^{r^{+}_i}, \quad i=1,2.
\end{equation}
Hence
\begin{equation} \label{Le6-14}    
\fm_{1,k_{1,1}} + m_{k_{1,1}} M_{1, \br_{k_{1,1}}}
  \equiv
 0    \mod  p_1^{r_{1,k_{1,1}}}, \; \;
 \fm_{2,k_{2,1}} + m_{k_{2,1}} M_{2, \br_{k_{2,1}}}
  \equiv
 0    \mod  p_2^{r_{2,k_{2,1}}}.
\end{equation}
According to \eqref{Le3-17}, \eqref{Le6-1} and \eqref{Le6-10}, we get that     $ \fm_{i,k_{i,j}} \in
I_{ p_i^{r_{i,k_{i,j}}-r_{i,k_{i,j-1}}}} $ and $r_{i,k_{i,1}} \\=r^{-}_{i}$, $k_{i,0}=r_{i,0} =0$.
Therefore
 $ \fm_{i,k_{i,1}} \in
I_{ p_i^{r_{i,k_{i,1}}}}$, $i=1,2$.

We fix $  m_{1}, m_2, \br_1,\br_2$.
Bearing in mind that $ I_M=[-[(M-1)/2],[M/2]] \cap \ZZ$ is a complete set of residues $\mod M$
 (see \eqref{Beg-1a}), we obtain that
there exists only one solution $ \fm_{i,k_{i,1}} \in
I_{ p_i^{r_{i,k_{i,1}}}}$
of  \eqref{Le6-14},   with $i \in \{1,2\}$. By  \eqref{Le6-10} and \eqref{Le6-12}, we have
$  p_2^{r_{2,k_{1,2}}}  \fm_{1,k_{1,2}}^{*}  =- m_{k_{1,2}}^{*}$,  $ \;\; p_1^{r_{1,k_{2,2}}} \fm_{2,k_{2,2}}^{*}  =- m_{k_{2,2}}^{*}$.
Substituting the above relations into  \eqref{Le7-0}, we obtain
\begin{multline}  \label{Le6-35}
\DD^{\ast,4}_{1, 1} \ll \sum_{\substack{k_{i,j}=1,2 \\i,j=1,2}} \;
\sum_{ (\br_1,\br_2) \in \UU_{ 1, 1 }} \; \sum_{\substack{ \beta_i \geq 0 \\ i=3,4}} \;
\sum_{\substack{  0< |\fm_{i,k_{i,2}}^{*}|  \leq n^{10} \\ i=1,2}}
 \;\prod_{i=1}^{2}\delta_{p_i^{r^{+}_i}}( \tilde{m}_{i} )
 \frac{ p_i^{-\beta_{i+2} -r_{i,k_{i^{'},2}}}}{ |\fm^{*}_{i,k_{i,2}}|^2 } \\
   \ll
 \sum_{\substack{r^{+}_i, r^{-}_i \in [1,n]\\ i=1,2}}
 \frac{1}{ p_1^{r^{-}_1} p_2^{r^{-}_2}}    \ll n^{2}.
\end{multline}
{\it Case} 5 : $\quad \lambda_{1} = 1$, $\lambda_{2}=0$.
 By \eqref{Le5-4} and \eqref{Le6-11},  we have
\begin{equation} \nonumber 
   r^{+}_{1} -r^{-}_{1} > \VV, \;\;   r^{+}_{2} -r^{-}_{2}    \leq \VV,
 \;\; r_{2,k_{1,2}}  \leq \VV \;\; \ad  \;\;
  \hat{r}_j = \max(r_{1,j}, r_{2,j}) > V>10 \VV, \; \;j=1,2.
\end{equation}
Hence $r^{+}_{2}=  \max(r_{2,1},r_{2,2})  \leq 2 \VV $.
Therefore  $  r_{1,j}=\hat{r}_j> V $  ($j=1,2$) and $  r^{-}_{1} = \min(r_{1,1},r_{1,2}) \geq V  $.
 By  \eqref{Le5-2}, we get
$  M_{1,\br_j}  p_2^{r_{2,j}} \equiv 1 \mod  p_1^{  r_{1,j}  }$  and
$  M_{1, \br_{k_{1,1}}}  p_2^{r_{2,k_{1,1}}} \equiv 1 \mod  p_1^{  r^{-}_1  }$,   $r^{-}_1 = r_{1, k_{1,1}}$.
In view of \eqref{Le6-13},
  we obtain $\fm_{1,k_{1,1}} + m_{k_{1,1}} M_{1, \br_{k_{1,1}}}  \equiv  0    \mod  p_1^{  r^{-}_1  }$
  and
\begin{equation} \nonumber \label{Le6-83}
       p_2^{r_{2,k_{1,1}}}  \fm_{1,k_{1,1}}  + m_{k_{1,1}}   \equiv
  0    \mod  p_1^{  r^{-}_1  }, \;  \ord_{p_1}(  p_2^{r_{2,k_{1,1}}}  \fm_{1,k_{1,1}}  + m_{k_{1,1}}  )  \geq  r^{-}_1 \geq  V.
\end{equation}
%
%
Let  $\beta_0 =\ord_{p_1}(m_{k_{1,1}})$,  $m_{k_{1,1}}^{*}=m_{k_{1,1}} p_{1}^{-\beta_0}$, $ (p_{1},m_{k_{1,1}}^{*} )  =1$ and  $\fm_{1,k_{1,1}}^{*} = \fm_{1,k_{1,1}} p_{1}^{-\beta_0}$.
We have  $\beta_0 \leq
  \log_{p_2} (\max(|m_1|, |m_2|)) \leq 10 \log_2 n  \leq  \VV/2$, and
   $r^{-}_1 -\VV/2 \geq V-\VV   \geq C_1 \log_2^2 n $ (see \eqref{Le5-4}).

Applying the Corollary, we get that the above congruence is equality :
\begin{equation} \label{Le6-75}    
    p_2^{r_{2,k_{1,1}}}  \fm_{1,k_{1,1}}^{*}  =- m_{k_{1,1}}^{*} \quad \ad \quad
    p_2^{r_{2,k_{1,1}}}  \fm_{1,k_{1,1}}  =- m_{k_{1,1}}, \quad \beta_0 =\ord_{p_1}(m_{k_{1,1}}).
\end{equation}
%
%
Now from \eqref{Le6-2}, we obtain, for $i=2$, that
\begin{equation} \label{Le6-51}    
   \fm_{2,k_{2,2}} + \fm_{2,k_{2,1}} p_2^{ r^{+}_2 -  r_{2,k_{2,1}}  } \equiv
A_0   \mod p_2^{r^{+}_2},
\end{equation}
with  $  A_0 = -m_{k_{2,2}} M_{2, \br_{k_{2,2}}} - m_{k_{2,1}} M_{2, \br_{k_{2,1}}}
 p_2^{ r^{+}_2 -  r_{2,k_{2,1}}  }  $.
We fix $  m_{1}, m_2, \br_1,\br_2$.
Bearing in mind that $ \fm_{2,k_{2,1}} \in
I_{ p_2^{r_{2, k_{2,1}}}}$  and $ \fm_{2,k_{2,2}} \in
I_{ p_2^{r^{+}_2 -r_{2,k_{2,1}}}} $ (see \eqref{Beg-1a} and \eqref{Le6-1}), we obtain that
there exists only one solution $(\fm_{2,1}, \fm_{2,2})  $
of  \eqref{Le6-51}.

Now substituting \eqref{Le6-10} and  \eqref{Le6-75} into  \eqref{Le7-0}, we obtain
\begin{multline}  \label{Le6-39}
\DD^{\ast,5}_{1, 0}= \sum_{\substack{k_{i,j}=1,2 \\i,j=1,2}} \;\;
\sum_{ (\br_1,\br_2) \in \UU_{ 1,  0 }} \;\;  \sum_{\substack{ \beta_{3i-3} \geq 0\\ i=1,2}}
\;\; \sum_{  0< |\fm_{1,k_{1,i}}^{*}|  \leq n^{10} , \; i=1,2}
 \;\; \prod_{i=1}^{2}\delta_{p_i^{r^{+}_i}}( \tilde{m}_{i} )
 \frac{  p_1^{-\beta_{3i-3}}     p_2^{ -r_{2,i}}}{ |\fm^{*}_{1,k_{1,i}}|^2 } \\  \ll
 \sum_{ r_{i,j} \in [1,n],\; i,j=1,2}
 \frac{1}{  p_2^{r_{2,1}+ r_{2,2}}} \ll n^{2}.
\end{multline}

From \eqref{Le4-1}, \eqref{Le7-0},  \eqref{Le6-35}  and \eqref{Le6-39},
 we obtain $\DD^{\ast}_{1, 1}=
   \DD^{\ast,4}_{1, 1} \ll n^2$ and
 $\DD^{\ast}_{1, 0} =\DD^{\ast,5}_{1, 0}
\ll n^2$.
The case $\quad \lambda_{1} = 0$, $\lambda_{2}=1$ is similar to the case $ \lambda_{1} = 1$, $\lambda_{2}=0$.
Therefore $\DD^{\ast}_{0, 1} \ll n^2$ and $\DD^{\ast}_{1, 1} + \DD^{\ast}_{1, 0} + \DD^{\ast}_{0, 1} \ll n^2$.
Hence Lemma 7 is proved. \qed

By  Lemma 6, Lemma 7 and \eqref{Le4-30},
the Theorem is proved. \qed \\

{\bf Remark 1.} The constant $C$ in the Theorem depends only on Yu's constant $C_1$ (see [Bu, p. 19]).\\

{\bf Remark 2.}  In [Le], we proved that Halton's sequence is of $L_q$-low discrepancy for all $s \geq 2$  and $q>0$. In [Le], we also proved the Central Limit Theorem  and moment convergence for
Hammersley's point set
\begin{equation}   \label{CLT1}
 \frac{ D(\bar{\bx}, \cH_{s+1,N} ) }{ \left\| D(\bar{\bx},\cH_{s+1,N}) \right\|_{2}}
	\stackrel{w}{\rightarrow} \cN(0,1), \qquad
 \frac{ D_{s+1,q}( \cH_{s+1,N} )}{ D_{s+1,2}( \cH_{s+1,N} )}
	\stackrel{N \rightarrow \infty}{\longrightarrow}   \kappa_q^{1/q}, \; s \geq 3,
\end{equation}
with $\kappa_q= \frac{1}{\sqrt{2 \pi}}\int_{-\infty}^{\infty} |u|^q e^{-u^2/2} d u$. In particular, we get that the lower bound \eqref{In7} is optimal for all $q>0$ for Hammersley's point set. The proof of these results is the same as the proof of the Theorem, but much more technical.
\\

{\bf Remark 3.}  The Halton sequence
is of $L_q$-low discrepancy for $s=1$ only for some sort of symmetrisation [KrPi]. It is very difficult to understand why  Halton's sequence is of $L_q$-low discrepancy  for $s\geq 2$ without any symmetrisation. The first idea is that the main tool (Theorem A) may be applied only for $s\geq 2$. But it is not possible that this explanation is complete. Because in [Le], we see the same problem for transition  from $s = 2$ to $s \geq 3$. Namely, we do not need any symmetrisation for the validity of  the Central Limit Theorem  for
Hammersley's point set for $s \geq 3$ (see \eqref{CLT1}). But we need a symmetrisation for the case $s = 2$. Thus, it remains to say that probably the auto-symmetrisation grows with the increase  of the dimension.

{\bf Address}: Department of Mathematics, Bar-Ilan University, Ramat-Gan, 5290002, Israel \\
{\bf E-mail}: mlevin@math.biu.ac.il\\
\end{document}